\documentclass[namedate,webpdf]{ima-authoring-template}

\journaltitle{}

\graphicspath{{Fig/}}

\theoremstyle{thmstyletwo}%
\newtheorem{theorem}{Theorem}
\newtheorem{proposition}[theorem]{Proposition}%
\newtheorem{corollary}[theorem]{Corollary}
\newtheorem{lemma}[theorem]{Lemma}

\newtheorem{example}[theorem]{Example}%
\newtheorem{remark}[theorem]{Remark}%

\numberwithin{equation}{section}

\def\C{\mathbb{C}}

\def\R{\mathbb{R}}

\def\span{\operatorname{span}}
\def\ker{\operatorname{ker}}
\def\vvec{\operatorname{vec}}

\usepackage{mathtools}

\DeclarePairedDelimiter{\norm}{\lVert}{\rVert}
\DeclarePairedDelimiter{\abs}{\lvert}{\rvert}

\newcommand{\subjectto}{\ensuremath{\ \text{s.t.}\ }}

\begin{document}

\DOI{10.48550/arXiv.2509.26344}
\copyrightyear{2026}
\pubyear{0000}
\access{Advance Access Publication Date: Day Month Year}
\appnotes{Paper}
\copyrightstatement{}


\title[Nearest matrix with multiple eigenvalues]{Nearest matrix with multiple eigenvalues by Riemannian optimization}

\author{Vanni Noferini\corresp[*]{Corresponding author: Vanni Noferini \href{email:vanni.noferini@aalto.fi}{vanni.noferini@aalto.fi}}\footnote{Supported by a Research Council of Finland grant (decision number 370932). Email: vanni.noferini@aalto.fi}\ORCID{0000-0002-1775-041X}
\address{\orgdiv{Department of Mathematics and Systems Analysis}, \orgname{Aalto University}, \orgaddress{\street{P.O. Box 11100}, \postcode{FI-00076}, \state{Aalto}, \country{Finland}}}}

\author{Lauri Nyman\footnote{Supported by an Engineering and Physical Sciences Research Council (EPSRC) grant EP/Z533786/1. Email:
lauri.nyman@manchester.ac.uk}\ORCID{0009-0002-9558-0951}
\address{\orgdiv{Department of Mathematics}, \orgname{University of Manchester}, \orgaddress{\postcode{M13 9PL}, \state{Manchester, England}, \country{UK}}}}
\author{Federico Poloni\footnote{Member of INDAM (Istituto Nazionale di Alta Matematica). Supported by ICSC--Centro Nazionale di Ricerca in High Performance Computing, Big Data, and Quantum Computing and by the PRIN project MOLE (2022ZK5ME7).
Email: federico.poloni@unipi.it}\ORCID{0000-0003-2291-990X}
\address{\orgdiv{Department of Computer Science}, \orgname{University of Pisa}, \orgaddress{\postcode{I-56127}, \state{Pisa}, \country{Italy}}}}

\authormark{Noferini, Nyman and Poloni}

\received{Date}{0}{Year}
\revised{Date}{0}{Year}
\accepted{Date}{0}{Year}
\makeatletter
\global\let\@history\@empty
\makeatother


\abstract{    Given a square complex matrix $A$, we tackle the problem of finding the nearest matrix with multiple eigenvalues or, equivalently when $A$ had distinct eigenvalues, the nearest defective matrix. To this goal, we extend the general framework described in [M. Gnazzo, V. Noferini, L. Nyman, F. Poloni, \emph{Riemann-Oracle: A general-purpose Riemannian optimizer to solve nearness problems in matrix theory}, Found. Comput. Math.] and based on variable projection and Riemannian optimization, allowing the ambient manifold to simultaneously track left and right eigenvectors. Our method also allows us to impose arbitrary complex-linear constraints on either the perturbation or the perturbed matrix; this can be useful to study structured eigenvalue condition numbers. We present numerical experiments, comparing with preexisting algorithms.}


\maketitle

\begin{center}
\textbf{To Nick}
\end{center}

\section{Introduction}

Matrix nearness problems consist in finding, given a matrix $A$ and a property $\mathfrak{P}$, a matrix $B$ that has property $\mathfrak{P}$ and is nearest to $A$, where distances are usually those induced by the Euclidean norm \cite{RO,nickthesis}. Nick Higham was a seminal figure in this field. Already his PhD thesis \cite{nickthesis} was devoted to this topic, and he continued to advance the knowledge on this line of research throughout his career \cite{nickspd,nicksurvey,nickcorr}.

In this paper, we address the following matrix nearness problem: For a square matrix $A \in \C^{n \times n}$, find the smallest (in the Frobenius norm) possible perturbation $\Delta \in \C^{n \times n}$ such that $A+\Delta$ has a multiple eigenvalue. This problem has a rich history. The distance from the set  of matrices with multiple eigenvalues is the same as the distance to the set of defective matrices; see Proposition \ref{prop:closure}. From a theoretical point of view, the problem of computing such distance is motivated by the study of eigenvalue condition numbers \cite{AlamBora05,ABBO,Demmel86, LE99, Wilk1, Wilk2}. Algorithmically, the problem of finding an efficient and reliable method to calculate the distance to defective matrices was posed in \cite{AlamBora05} and further studied by various authors \cite{AFS,ABBO,BGMN,gugbook,Russians23}. Even though, in \cite{AFS,AlamBora05, ABBO,BGMN,gugbook,Russians23}, the emphasis is on defective matrices, we prefer to (equivalently) focus on matrices with multiple eigenvalues because it is always possible to find a nearest one to every input $A$. This is in contrast with the fact that there exist some $A$ that lack a nearest defective matrix. On the other hand, all these problematic inputs correspond to instances where the problem is in some sense trivial, because the sought distance is $0$ and the sought perturbation is $\Delta=0$, but at the same time $A$ is not defective, and hence no nearest defective matrix to it exists; see \cite[Theorem 6]{ABBO} and Proposition \ref{prop:fullchar}. Hence, at least assuming that in practice one is interested in inputs for which the solution exists, our analysis and our algorithm are also applicable to the problem of finding a nearest defective matrix.

For this problem, we propose a method based on Riemannian optimization and variable projection. More specifically, we adapt the general framework described in \cite{RO}, thus translating the problem of computing the nearest matrix with multiple eigenvalues into a minimization task on the complex Stiefel manifold $V_2(\C^n)$, i.e., the manifold of $n \times 2$ complex matrices with orthonormal columns. An important novelty with respect to most previously existing methods is that our approach allows us to even add an additional linear constraint on either $\Delta$ or $A+\Delta$ (or both). In other words, we can require that the perturbation and/or the perturbed matrix belongs to a prescribed complex subspace of $\C^{n \times n}$. For example, we can find the smallest Toeplitz perturbation to obtain a matrix with a double eigenvalue, or the nearest matrix with a double eigenvalue and a certain prescribed sparsity pattern, and so on. This feature is applicable, for example, to the study of \emph{structured} eigenvalue condition numbers. Another consequence of this flexibility is that our algorithm can also be applied, via companion matrices, to the problem of computing the nearest polynomial with a double root. To our knowledge, the only previous work that considered additional constraints on the problem --- though in a different sense than we do here --- was \cite{BGMN,gugbook}, where real perturbations were also studied. While the approach described  in~\cite[Section~VI.3]{gugbook} can in principle also deal with more general structures, in practice only the case of the real structure is considered there.

The structure of the paper is as follows.  In Section \ref{sec:problem}, we formally describe and analyze the problem, and we describe a Riemannian optimization method to solve it, including the case where additional linear constraints are considered. In Section \ref{sec:2}, we move our attention back to the unstructured (that is, no additional constraint) version of the problem, as formulae simplify in that case. In Section \ref{sec:ne} we present numerical experiments, and we finally draw some conclusions in Section \ref{sec:conclusions}.

\section{Statement of the problem and reformulation as optimization on manifolds}\label{sec:problem}

\subsection{Nearest matrix with multiple eigenvalues and nearest defective matrix}\label{sub:rant}

Recall that a square matrix is defective if it possesses at least one eigenvalue whose algebraic multiplicity is strictly larger than its geometric multiplicity. By definition, it is clear that a defective matrix has multiple eigenvalues, but the reverse implication is false. For example, $0$ is a multiple eigenvalue of the zero matrix (of size at least $2$), but the zero matrix is not defective. Having fixed the size $n \geq 2$, denote below by $\mathcal{M} \subset \C^{n \times n}$ the set of matrices with multiple eigenvalues and by $\mathcal{D} \subset \C^{n \times n}$ the set of defective matrices. The previous remarks, together with the observation that some (in fact, almost all) matrices do not have multiple eigenvalues, yield the strict inclusions $\mathcal{D} \subsetneq  \mathcal{M} \subsetneq \C^{n \times n}.$

In some contexts, for example the analysis of eigenvalue condition numbers, it is of interest to compute the distance from a matrix to the set $\mathcal{D}$. One may also want to find a matrix in $\mathcal{D}$ that minimizes the distance to the input. Problematically, however, $\mathcal{D}$ is not closed. When $n=2$, this can be easily seen by considering the converging sequence
\[  (A_k)_k \subset \mathcal{D}, \quad (A_k) = \begin{bmatrix}
    0 & \frac{1}{k+1}\\
    0 & 0
\end{bmatrix}\]
whose limit does not belong to $\mathcal{D}$. It is clear that similar examples also exist for higher sizes $n \geq 3$. Following \cite[Section 2]{NN25}, we observe that in a metric space equipped with a distance $d$ it is always possible to define the distance from an element $x$ to a set $\mathcal{S}$, as the infimum over all $y \in \mathcal{S}$ of $d(x,y)$, i.e., the distance between $x$ and $y$. More concretely, in the present article our metric space is $\C^{n \times n}$ and our distance is $d(X,Y)=\|X-Y\|_F$. However, the existence of at least one element of $\mathcal{S}$ where the infimum is achieved is only guaranteed if $\mathcal{S}$ is closed; moreover, for all $x$ and all $\mathcal{S}$, one has that $d(x,\mathcal{S})=d(x,\overline{\mathcal{S}}$). We therefore conclude that, while the problem of finding the distance to the set $\mathcal{D}$ is well defined, computing a nearest defective matrix to a given one is, in general, ill-posed. Fortunately, the cure is very simple, because one can look instead for a nearest element in $\overline{\mathcal{D}}$.

\begin{proposition}\label{prop:closure}
    $\overline{\mathcal{D}}=\mathcal{M}$.
\end{proposition}
\begin{proof}
By continuity of the characteristic polynomial and the discriminant, it is clear that $\overline{\mathcal{D}} \subseteq \mathcal{M}$. For the reverse inclusion, take $A \in \mathcal{M} \setminus \mathcal{D}$ and let $\lambda_0$ be a multiple eigenvalue of $A$. As $A \notin \mathcal{D}$, there exists a similarity $X \in GL(n,\C)$, and a matrix $B \in \C^{(n-2) \times (n-2)}$ such that $A=X(\lambda_0 I_2 \oplus B)X^{-1}$. It then suffices to consider $(A_k)_k \subset \mathcal{D}$ with $$A_k = X\left(\begin{bmatrix}
    \lambda_0 & \frac{1}{k+1}\\
    0 & \lambda_0
\end{bmatrix} \oplus B\right)X^{-1} \Rightarrow A_k \rightarrow A.$$
\end{proof}

Note that it is easy to construct concrete examples of a matrix whose nearest defective matrix does not exist. Take an arbitrary matrix $A \in \mathcal{M} \setminus \mathcal{D}$ (for example, the zero matrix $A=0$). Then, clearly $0=d(0,\mathcal{D})=\inf_{Y \in \mathcal{D}} \|Y\|_F$, but the infimum cannot be achieved at any element of $\mathcal{D}$ because $Y \in \mathcal{D} \Rightarrow Y \neq 0 \Rightarrow \|Y\|_F > 0$. Indeed, we argue below that these are the only counterexamples. This fact is essentially a corollary of \cite[Theorem 6]{ABBO}.

\begin{proposition}\label{prop:fullchar}
    There exists a matrix in $\mathcal{D}$ nearest to $A$ if and only if $A \not\in \mathcal{M} \setminus \mathcal{D}$.
\end{proposition}
\begin{proof}
    If $A \not \in \mathcal{M}$ then the statement follows by \cite[Theorem 6]{ABBO}. If $A \in \mathcal{M}$, then by Proposition \ref{prop:closure} it holds $d(A,\mathcal{D})=d(A,\mathcal{M})=0$. Over $\mathcal{M}$, the unique minimizer of the distance is $A$; thus, a minimizer does not exist in $\mathcal{D}$ if and only if $A \in \mathcal{M}\setminus \mathcal{D}$.
\end{proof}

In summary, what this digression shows is that the problem of computing the nearest defective matrix is completely equivalent to that of computing the nearest matrix with multiple eigenvalues for \emph{almost every input}, in the technical sense that the exceptional set of inputs for which a nearest defective matrix does not exist has zero measure. Arguably, the equivalence of these problems holds for every input of practical interest, because all the exceptions belong to $\mathcal{M}$, and hence have distance $0$ both from $\mathcal{M}$ and from $\mathcal{D}$. However, to avoid bothering with special cases in our statements, in the following we focus only on the problem of finding a nearest matrix in $\mathcal{M}$.

\subsection{The nearest matrix with multiple eigenvalues as an optimization problem}

We now turn to describing the set $\mathcal{M}$. The quality of having multiple eigenvalues can be characterized in terms of the eigenvectors of the matrix. This result is well known; see, e.g., \cite[Theorem 1.4.12]{hornjohnson}. We state a slightly different version in Lemma \ref{lemma:multiple_eval}; up to minor tweaks, our proof is essentially the same as in \cite[Lemma 1]{Russians23}, but we include it for completeness. Here and throughout, we denote by $K^*$ the conjugate transpose of the matrix (or vector) $K$, by $K^T$ its pure transpose (without conjugations), and by $\overline{K}$ its elementwise conjugate.

\begin{lemma}\label{lemma:multiple_eval}
    Let $A \in \C^{n \times n}$ and $\lambda \in \C$. Then, $\lambda$ is a multiple eigenvalue of $A$ if and only if there exist orthonormal vectors $u,v \in \C^n$ such that $u^* A = \lambda u^*$ and $A v=\lambda v$.
\end{lemma}
\begin{proof}
    Suppose first $\lambda$ is an eigenvalue of $A$ with algebraic multiplicity $\geq 2$. Then, there exists a Schur form of $A$, say, $T=QAQ^*$, such that $e_1^T T e_1 = \lambda = e_n^T T e_n$. Hence, $T e_1 = \lambda e_1$ and $e_n^T T = \lambda e_n^T$, and it suffices to define $v:=Q^* e_1$ and $u^*:= e_n^T Q$ to obtain vectors with the desired properties.

Conversely, assume that $u,v$ exist as in the statement; this immediately implies that $\lambda$ is an eigenvalue.  Suppose for a contradiction that $\lambda$ has algebraic multiplicity $1$. Then, its geometric multiplicity is also $1$. Thus, $\span(u)=\ker((A-\lambda I)^*)$ and hence $v \in \span(u)^\perp=\mathrm{colspace}(A-\lambda I)$, implying in turn the existence of $w$ such that $v=Aw-\lambda w$. It follows that $(v,w)$ is a Jordan chain at $\lambda$ for $A$, contradicting the assumption that $\lambda$ is simple.
\end{proof}

In addition, the following characterization was proved in~\cite[Lemma~1, Theorem~5, Theorem~6]{ABBO}.
\begin{theorem} \cite{ABBO}\label{thm:svdcharacterization}
    Let $A \in \C^{n \times n}\setminus \mathcal{M}$ be a matrix without multiple eigenvalues, and let $\delta=d(A,\mathcal{M}) > 0$ (where the distance is measured in the Frobenius norm). Then there exists $\lambda \in \C$ such that, for a suitable choice of a singular value decomposition $A-\lambda I = \sum_{i=1}^n u_i \sigma_i v_i^*$, it holds for $\Delta=-\delta u_n v_n^*$ that:
    \begin{itemize}
        \item[(i)] $A+\Delta$ is a matrix in $\mathcal{D} \subset \mathcal{M}$ nearest to $A$;
        \item[(ii)] $\lambda$ is a multiple eigenvalue of $A+\Delta$;
        \item[(iii)] $\delta = \sigma_n$;
        \item[(iv)] $u_n^*v_n=0$;
        \item[(v)] $(A+\Delta-\lambda I)v_n = (A+\Delta-\lambda I)^*u_n=0$.
    \end{itemize}
     \end{theorem}
Note that, in the statement of Theorem \ref{thm:svdcharacterization}, the choice of the vectors $u_n$ and $v_n$ in the SVD is not unique if $\sigma_{n-1}=\sigma_n$; in this scenario, Theorem \ref{thm:svdcharacterization} does not claim that all possible choices satisfy the statement, but rather that one ``good" choice exists.

Based on Lemma \ref{lemma:multiple_eval}, we can reformulate the problem of finding the nearest matrix with multiple eigenvalues as an optimization problem over an orthonormal pair $u,v \in \C^n$. To this end, let $V_2(\C^n)$ denote the complex Stiefel manifold of orthogonal 2-frames in $\C^n$. In the following, we refer to this manifold more succinctly as $V_2$. We can then express the problem of finding the nearest matrix with multiple eigenvalues as the nested minimization problem
\begin{equation} \label{problemwithouttransposes}
    \min_{[u \ v] \in V_2} \min_{\lambda \in \C} \min_{{\Delta\in\mathcal{S}}} \, \norm{\Delta}_F^2 \subjectto (A+\Delta-\lambda I)v = 0, \, (A+\Delta-\lambda I)^* u = 0.
\end{equation}
In \eqref{problemwithouttransposes}, the set $\mathcal{S}$ denotes a desired complex linear subspace of $\C^{n \times n}$.

We make immediately a small but important change in our formulation of~\eqref{problemwithouttransposes}: We take the elementwise conjugate of the second constraint, and define
\begin{equation} \label{fuvl}
    f(u,v,\lambda) = \min_{{\Delta\in\mathcal{S}}} \, \norm{\Delta}_F^2 \subjectto \begin{bmatrix} (A+\Delta-\lambda I)v \\ (A+\Delta-\lambda I)^T \overline{u} \end{bmatrix} = 0.
\end{equation}

Clearly, this modification does not alter the problem; however, it has the advantage that now the constraint is linear in the entries of $\Delta$. Linearity will be instrumental in Subsection~\ref{sec:closedform}, where we exploit it to write the solution of the minimum norm problem~\eqref{fuvl} explicitly, using a pseudoinverse. The problem \eqref{problemwithouttransposes} now becomes $$\min_{[u \ v] \in V_2} \min_{\lambda \in \C} f(u,v,\lambda).$$ 

As analyzed in \cite[Section 2]{RO}, functions that define the solution of a minimum-norm problem like~\eqref{fuvl} may be discontinuous. This discontinuity can be circumvented with a regularization technique, again following~\cite{RO}. The Tikhonov regularized least-squares solution corresponds to the minimizer of 
\begin{equation} \label{fuvlreg}
f_\varepsilon(u,v,\lambda) = \min_{{\Delta\in\mathcal{S}}} \| \Delta \|_F^2 + \varepsilon^{-1}  \left\| \begin{bmatrix} (A+\Delta-\lambda I)v \\ (A+\Delta-\lambda I)^T \overline{u} \end{bmatrix}
  \right\|^2 .
\end{equation}
Note that we do not relax the constraint $\Delta \in \mathcal{S}$, which always holds exactly. 

We solve the problem 
\begin{equation} \label{minfuvepsilon}
    \min_{[u \ v] \in V_2} \min_{\lambda \in \C} f_\varepsilon (u,v,\lambda)
\end{equation}
multiple times, with a decreasing sequence of parameters $\varepsilon$, using the minimizer of the previous iteration as a starting point for the next iteration. Under suitable assumptions, the solutions of this sequence of regularized minimum problems converge to a solution of the non-regularized problem \eqref{problemwithouttransposes} in the limit $ \varepsilon \rightarrow 0$. This algorithm is an instance of the \emph{penalty method} from constrained optimization; its convergence properties are described in more detail in \cite[Section 4]{RO}. The authors of \cite{RO} also discussed another solution algorithm with stronger convergence properties, the \emph{augmented Lagrangian method}. For the problem studied in this paper, in the augmented Lagrangian method we have at each step $\varepsilon>0$ and a vector of Lagrange multipliers $y\in\mathbb{C}^{2n}$, and we need to compute
\begin{equation} \label{augmented_lagrangian_method}
\min_{[u,v]\in V_2} \min_{\lambda\in\mathbb{C}}  f_{\varepsilon,y}(u,v,\lambda), \quad
f_{\varepsilon,y}(u,v,\lambda) = \norm{\Delta}_F^2 + \varepsilon^{-1} \left\| \begin{bmatrix} (A+\Delta-\lambda I)v \\ (A+\Delta-\lambda I)^T \overline{u} \end{bmatrix} +y\varepsilon
  \right\|^2.
\end{equation}
We refer the reader once again to~\cite[Sections~2.3 and~4]{RO} for details on outer iteration steps of the augmented Lagrangian method, and we focus in the next subsection on giving a closed form for the solution of~\eqref{augmented_lagrangian_method}. Note that the Tikhonov regularization approach \eqref{fuvlreg} can be seen as a special case of the augmented Lagrangian approach \eqref{augmented_lagrangian_method}, obtained by setting $y=0$ at every iteration.

To conclude this subsection, following \cite[Remark 2.6]{RO} we observe that the machinery described above for linear constraints of the form $\Delta \in \mathcal{S}$ is also applicable when we require instead that $A+\Delta$ belongs to a subspace $\mathcal{S}$. Indeed, note the unique decomposition $A=A_\mathcal{S}+A_\perp$ with $A_\mathcal{S} \in \mathcal{S}$ and $A_\perp \in \mathcal{S}^\perp$. Suppose now $A+\Delta \in \mathcal{S}$. Then, defining $\Delta_\mathcal{S}=A_\perp+\Delta \in \mathcal{S}$, we have that $\|\Delta\|_F^2 = \|\Delta_\mathcal{S}\|^2_F + \| A_\perp\|_F^2$, and therefore minimizing the norm of $\Delta$ such that $A+\Delta \in \mathcal{S} \cap \mathcal{M}$ is equivalent to minimizing the norm of $\Delta_\mathcal{S} \in \mathcal{S}$ such that $A_\mathcal{S}+\Delta_\mathcal{S} \in \mathcal{M}$.

\subsection{Closed form for the objective function} \label{sec:closedform}
In this subsection, we wish to give a closed form for the objective function of the regularized problem  $f_{\varepsilon,y}(u,v,\lambda)$, following the technique (sometimes called \emph{variable projection}) given in \cite{UseM14} and \cite[Theorem 2.4]{RO} for similar problems.

Consider the Frobenius inner product on $\C^{n \times n}$, defined as $\langle X,Y \rangle = \mathrm{trace}(Y^*X)$, and note that $\|X\|_F=\sqrt{\langle X,X\rangle}$. Let $P^{(1)}, \dots, P^{(p)} \in \mathbb{C}^{n\times n}$ be an orthonormal basis for the subspace $\mathcal{S}$, where $p = \dim \mathcal{S}$. We can then write
\[
\Delta = \sum P^{(j)}\delta_j, \quad \vvec \Delta = \mathcal{P}\delta, \quad \mathcal{P} = \begin{bmatrix}
    \vvec P^{(1)} & \vvec P^{(2)} & \dots & \vvec P^{(p)}
\end{bmatrix} \in \mathbb{C}^{n^2 \times p},
\]
and orthonormality guarantees that $\mathcal{P}^T\mathcal{P}=I_p$ and $\norm{\Delta}_F = \norm{\delta}$.

We start with a useful lemma.
\begin{lemma}
    Let $u,v\in\mathbb{C}^n$. The unique matrix $M(u,v) \in\mathbb{C}^{2n\times p}$ such that
    \[
    M(u,v)\delta = \begin{bmatrix}
        \Delta v\\
        \Delta^T \overline{u}
    \end{bmatrix}
    \]
    for every $\Delta\in\mathcal{S}$ is given by
    \begin{equation} \label{defM}
    M(u,v) = 
    \begin{bmatrix}
        P^{(1)}v & P^{(2)}v & \dots & P^{(p)}v \\
        {P^{(1)}}^T \bar u & {P^{(2)}}^T \bar u & \dots & {P^{(p)}}^T \bar u
    \end{bmatrix}
    =
    \begin{bmatrix}
        v^T \otimes I_n\\
        I_n\otimes u^*
    \end{bmatrix}\mathcal{P}.        
    \end{equation}
\end{lemma}
\begin{proof}
    The first expression for $M(u,v)$ follows from $\Delta = \sum P^{(j)}\delta_j$ and $\Delta^T = \sum {P^{(j)}}^T\delta_j$. The second one is obtained noting that for every matrix $P^{(j)}$ we have
    \[
    P^{(j)}v = \vvec P^{(j)}v = (v^T\otimes I_n) \vvec P^{(j)}, \quad
    {P^{(j)}}^T \overline{u} = \vvec u^*P^{(j)} = (I_n\otimes u^*) \vvec P^{(j)}
    \]
    by the properties of Kronecker products.
\end{proof}

\begin{theorem} \label{thm:fepsilonv}
Let $u,v\in \mathbb{C}^n$, $y \in \C^{2n}$, and $\varepsilon>0$ be given. Let $M=M(u,v)$ be as in~\eqref{defM}, and
\begin{equation} \label{defr}
    r = r(u,v)= \begin{bmatrix}
      (\lambda I - A)v\\
      (\lambda I - A)^T\overline{u}
  \end{bmatrix} - y\varepsilon.
\end{equation}
Then,
\begin{equation} \label{eq:feps}
f_{\varepsilon,y}(u,v,\lambda) = r^*(MM^* + \varepsilon I)^{-1}r.    
\end{equation}
Moreover, the matrix that gives the minimum in~\eqref{augmented_lagrangian_method} is $\Delta_* = \sum_{i=1}^p P^{(i)} (\delta_*)_i$, with
\begin{equation} \label{deltastar}
\delta_* = M^*(MM^*+\varepsilon I_{2n})^{-1}r = (M^*M+\varepsilon I_{p})^{-1}M^*r.
\end{equation}
\end{theorem}
\begin{proof} 
We have
\begin{equation} \label{Mdeltaminusr}
\begin{bmatrix} (A+\Delta-\lambda I)v\\ (A+\Delta-\lambda I)^T \overline{u} \end{bmatrix}+y\varepsilon = \begin{bmatrix} \Delta v \\ \Delta^T \overline{u} \end{bmatrix}  + \begin{bmatrix} (A-\lambda I)v \\ (A-\lambda I)^T \bar u \end{bmatrix} +y\varepsilon = M\delta-r,
\end{equation}
using the definitions in \eqref{defM} and \eqref{defr}. Given that $\norm{\Delta}_F = \norm{\delta}$, we can express the objective function in~\eqref{augmented_lagrangian_method} as
\[
f_{\varepsilon,y}(u,v,\lambda) = \min_{{\delta\in\mathbb{C}^p}} \; \norm{\delta}^2 + \varepsilon^{-1}\norm{M\delta-r}^2.
\]
Defining $\delta_*$ as in~\eqref{deltastar}, a tedious but straightforward manipulation yields the identity
\[
\norm{\delta}^2 + \varepsilon^{-1}\norm{M\delta-r}^2 = \varepsilon^{-1}(\delta-\delta_*)^*(M^*M+\varepsilon I)(\delta-\delta_*) + r^*(MM^*+\varepsilon I)^{-1}r.
\]
Since $M^* M+\varepsilon I$ is positive definite, this formula shows that the objective function has minimum $r^*(MM^*+\varepsilon I)^{-1}r$, achieved when $\delta =\delta_*$.
\end{proof}

We note that Theorem \ref{thm:fepsilonv} can be applied also to find closed form and minimizers for~\eqref{fuvlreg}, by setting $y=0$, and for~\eqref{fuvl}, by additionally taking the limit $\varepsilon\to 0$.
\begin{corollary}
    Using the same notation~\eqref{defM} and~\eqref{defr}, the minimum of~\eqref{fuvl} is computed over a non-empty feasible set if and only if $\lim_{y\to0,\varepsilon \to 0}r\in\operatorname{Im} M$, and in this case
    \[
    f(u,v,\lambda) = \norm{M^+r}^2 = r^*(MM^*)^+r = \lim_{y\to0,\varepsilon\to 0} r^*(MM^*+\varepsilon I)^{-1}r.
    \]
\end{corollary}
\begin{proof}
    Using~\eqref{Mdeltaminusr}, we reduce the minimization to a minimum-norm problem
    \[
    \min \norm{\delta}^2 \text{ s.t. } M\delta = r,
    \]
    for which these results are standard. See, e.g., \cite[Theorem~1.2.10]{Bjorck}.
\end{proof}

\begin{remark}
Unfortunately, it seems that one cannot easily express the computations in this section only with the symbol $\cdot^*$, avoiding unpaired transpositions $\cdot^T$ and conjugations $\overline{\cdot}$. Indeed, if we do not include the transposition in the second block of~\eqref{fuvlreg}, then the analogue of~\eqref{Mdeltaminusr} reads
\begin{equation} \label{otherversion}
\begin{bmatrix}
(A+\Delta-\lambda I)v\\
(A+\Delta-\lambda I)^*u
\end{bmatrix}
=
\begin{bmatrix}
M_1\\ 0
\end{bmatrix}\delta +
\begin{bmatrix}
0\\M_2
\end{bmatrix}\overline{\delta} - r.
\end{equation}
At this point we no longer have a least-squares problem with a linear constraint: the constraint in~\eqref{otherversion} is linear in $\operatorname{Re}(\delta)$ and $\operatorname{Im}(\delta)$ separately, and we need to separate real and imaginary parts to use the standard results on least-squares problems with pseudoinverses.
\end{remark}

The following result holds for the minimizer.
\begin{lemma} \label{lem:minimizer_Prank2}
    The minimizer $\Delta_*$ is the orthogonal projection onto $\mathcal{S}$ (in the Frobenius inner product) of a matrix with rank at most 2.
\end{lemma}
\begin{proof}
Let us set
\begin{equation} \label{defz}
z = \begin{bmatrix}
    z_v\\
    \overline{z_u}
\end{bmatrix} = (MM^*+\varepsilon I)^{-1}r.    
\end{equation}
Then,
\[
\delta_* = M^* z = \mathcal{P}^* \begin{bmatrix}
    \overline{v} \otimes I_n & I_n \otimes u
\end{bmatrix} \begin{bmatrix}
    z_v\\
    \overline{z_u}
\end{bmatrix},
\]
and hence
\begin{equation} \label{rank2prop}
\vvec \Delta_* = \mathcal{P}\delta_* = \mathcal{P}\mathcal{P}^* \left(\overline{v}\otimes z_v + \overline{z_u} \otimes u \right) = \mathcal{P}\mathcal{P}^* \vvec \left(z_v v^* + u z_u^* \right).    
\end{equation}
Since $\mathcal{P}$ has orthonormal columns, the matrix $\mathcal{P}\mathcal{P}^*$ is the orthogonal projector onto $\mathcal{S}$.
\end{proof}

We note in passing that, based on Theorem \ref{thm:svdcharacterization}, in the limit $y=0$ and $\varepsilon \rightarrow 0$, and at least at the global minimum and in the unstructured case where $\mathcal{S}=\C^{n \times n}$, we expect that the minimizer has rank $1<2$, i.e., a stronger rank property than predicted by Lemma \ref{lem:minimizer_Prank2}. Indeed, and remarkably, this observation actually also holds for $\varepsilon > 0$, as we shall prove later in Theorem \ref{thm:specialgradient}.

It turns out that there exists a closed form solution formula for the minimizer of the function \eqref{eq:feps} over $\lambda$. We give this formula in Theorem \ref{thm:minimizerlambda}.

\begin{theorem} \label{thm:minimizerlambda}
    Let $M,r$ be defined as in \eqref{defM} and \eqref{defr}, and let $f_{\varepsilon,y}(u,v,\lambda)$ be defined as in~\eqref{augmented_lagrangian_method}. Define  
    \begin{equation}
       r_1 =  \begin{bmatrix}
      v\\
      \overline{u}
  \end{bmatrix}, \quad r_0 = -\begin{bmatrix}
      A v \\
      A^T \overline{u}
  \end{bmatrix} - y\varepsilon,
    \end{equation} so that $r=r_1\lambda + r_0$, and
    \begin{equation} \label{eq:Nabc}
        a= r_1^*(MM^*+\varepsilon I)^{-1}r_1, \ b=r_1^* (MM^*+\varepsilon I)^{-1} r_0, \ c=r_0^*(MM^*+\varepsilon I)^{-1}r_0.
    \end{equation} 
    Then, the minimizer $\lambda_* = \arg \min_{\lambda \in \C} f_{\varepsilon,y}(u,v,\lambda)$ is $\lambda_* =-\frac{b}{a}$, with minimum $f_\varepsilon(u,v,\lambda_*) = \frac{ac-|b|^2}{a}$. 
\end{theorem}
\begin{proof}
We can express \eqref{eq:feps} as
\begin{equation} \label{f_as_quadratic_in_lambda}
    f_{\varepsilon,y}(u,v,\lambda) = \begin{bmatrix} \lambda^* & 1 \end{bmatrix} \begin{bmatrix}
        a & b\\
        \overline{b} & c
    \end{bmatrix}   \begin{bmatrix}
        \lambda\\
        1
    \end{bmatrix},
\end{equation} 
and it is easy to verify that the minimizer of this quadratic function is $\lambda_*=-\frac{b}{a}$, with minimum $\frac{ac-|b|^2}{a}$.
\end{proof}

Theorem \ref{thm:minimizerlambda} shows that the objective function becomes \[f_{\varepsilon,y}(u,v) = \min_{\lambda \in \C} r^*(MM^*+\varepsilon I)^{-1} r= \frac{ac-|b|^2}{a},\]
where $a,b,c$ are defined as in \eqref{eq:Nabc}. In this form, the problem becomes a Riemannian optimization problem over the manifold $V_2(\C^n)$. 

In Theorem \ref{thm:gradient}, we calculate the Euclidean gradient of the objective function $f_{\varepsilon,y}(u,v)$. The Riemannian gradient can then be computed by projecting the Euclidean gradient onto the tangent space of the Riemannian manifold in question \cite[Section~3]{Boumal}, in this case $V_2(\C^n)$.
\begin{theorem} \label{thm:gradient}
    Let $M = M(u,v)$ be as in \eqref{defM} and \eqref{defr}, let the minimizer $\lambda_*$ be as in Theorem \ref{thm:minimizerlambda}, let $z_u,z_v$ be as in~\eqref{defz}, and let the minimizers $\delta_*,\Delta_*$  be as in~\eqref{rank2prop}.
    Then, the Euclidean gradient of $f_{\varepsilon,y}(u,v) = \min_{\lambda \in \C} r^*(MM^*+\varepsilon I)^{-1}r$ is
    \begin{equation} \label{gradf}
         \nabla f_{\varepsilon,y}\left(\begin{bmatrix}
    u & v
\end{bmatrix}\right) = 2\begin{bmatrix}
    (\lambda_* I - A - \Delta_*)z_u &
    (\lambda_* I - A - \Delta_*)^* z_v
\end{bmatrix}.
    \end{equation}
\end{theorem}
\begin{proof}
To compute the gradient, we can mimic the proof of~\cite[Appendix~A]{RO}: when $v,u$ are modified to $v+t\dot{v}, u + t\dot{u}$, we have
\[
\dot{r} = \begin{bmatrix}
    (\lambda_* I -A)\dot{v}\\
    (\lambda_* I - A^T)\overline{\dot{u}}
\end{bmatrix}, \quad \dot{M}\delta_* = 
\begin{bmatrix}
    \vvec(\Delta_* \dot{v}) \\
    \vvec(u^*\Delta_*)
\end{bmatrix} = 
\begin{bmatrix}
    \vvec(\Delta_* \dot{v}) \\
    \vvec(\Delta_*^T \overline{\dot{u}})
\end{bmatrix}.
\]
Hence
\[
\dot{f}_{\varepsilon,y} = 2\operatorname{Re} z^* \begin{bmatrix}
    (\lambda_* I - A - \Delta_*)\dot{v}\\
    (\lambda_* I - A^T - \Delta_*^T)\overline{\dot{u}}
\end{bmatrix}
\]
and
\[
\nabla f_{\varepsilon,y}\left(\begin{bmatrix}
    u & v
\end{bmatrix}\right) = 2\begin{bmatrix}
    (\lambda_* I - A - \Delta_*)z_u &
    (\lambda_* I - A - \Delta_*)^* z_v
\end{bmatrix}.
\]
Note that in this computation we can consider $\lambda_*$ as a constant: indeed, if $\lambda_*$ is the minimizer of~\eqref{f_as_quadratic_in_lambda}, then because of this minimality we have
\[
\frac{df(\lambda,u,v)}{d\lambda}\bigg|_{\lambda_*} = 0,
\]
and hence the contribution of $\lambda_*$ in the first-order change to $f$ vanishes.    
\end{proof}

\subsection{Two examples}

\begin{example}
    Let $A=\begin{bmatrix}
        1&0\\
        0&0
    \end{bmatrix}$. By fixing otherwise redundant phases, we may assume that $e_1^Tu \in \R, e_1^Tv \in \R$. Hence, without loss of generality we can parametrize
    \[ u=\begin{bmatrix}
        \cos(t)\\
        -e^{i \theta} \sin(t)
    \end{bmatrix}, \quad v=\begin{bmatrix}
        \sin(t)\\
        e^{i \theta} \cos(t)
    \end{bmatrix}.    \]
Then,
\[  M = \begin{bmatrix}
\sin(t) & 0 & e^{-i \theta} \cos(t) & 0\\
0&\sin(t)&0&e^{-i \theta} \cos(t)\\
\cos(t)&-e^{i \theta} \sin(t)&0&0\\
0&0&\cos(t)&-e^{i \theta} \sin(t)
\end{bmatrix}, \qquad r=\begin{bmatrix}
     (\lambda-1) \sin(t)\\
   \lambda e^{i \theta} \cos(t)\\
    (\lambda-1)\cos(t)\\
-\lambda e^{-i \theta}\sin(t)
\end{bmatrix}\]
Some tedious, but straightforward, algebraic manipulations lead us to
\[  a=\frac{2}{1+\varepsilon}, \quad b=-\frac{1}{1+\varepsilon}, \quad c=\frac{7+\cos(4t)+4\varepsilon}{8+12\varepsilon+4\varepsilon^2} \Rightarrow -\frac{b}{a} = \frac12 \in \R.\]
Substituting this optimal value for $\lambda$, and taking the limit $\varepsilon \rightarrow 0$, yields the function (independent of $\theta$!)
\[  f(u,v)=\frac{\cos(4t)+3}{8} \]
which is minimized when $\cos(4t)=-1 \Rightarrow t=\frac{\pi}{4}+\frac{k \pi}{2}$.

Computing the optimal $\Delta=\lim_{\varepsilon \rightarrow 0}\vvec^{-1}(M^*(MM^*+\varepsilon I)^{-1}r)$, where in $r$ and in $M$ we have set $\lambda=\frac12$ and $t \in \{ \frac{\pi}{4}, \frac{3\pi}{4}, \frac{5\pi}{4}, \frac{7\pi}{4}  \}$, yields the family of optimal perturbations

\[ \Delta = \frac14 \begin{bmatrix} -1 & -e^{-i \theta} \\
e^{i \theta} & 1
\end{bmatrix} \quad \theta \in [0,2 \pi[ . \]
Note that these optimal perturbations $\Delta$ are all rank $1$, and that this family of complex solutions also contains two real solutions for $\theta \in \{0,\pi \}$. Moreover, all the resulting perturbed matrices $A+\Delta$ are similar to each other via a diagonal and unitary similarity matrix, and have a double eigenvalue $\frac12$.

\end{example}

\begin{example} \label{example:companion}
    In this example, we seek the monic scalar polynomial of degree $2$ having a double root and nearest to $z^2-z$. Here, given two degree-$k$ monic polynomials $a(z)=z^k+\sum_{i=0}^{k-1} a_i z^i$ and $b(z)=z^k+\sum_{i=0}^{k-1}b_i z^i$, we define their squared distance to be $\sum_{i=0}^{k-1} |a_i-b_i|^2$. Of course, this toy example can also be solved directly by looking for the nearest polynomial of the form $(z-\eta)^2$ and optimizing over the scalar complex parameter $\eta$. However, since it can be solved analytically even with our method, let us anyway analyze it for illustration purposes.

    The problem is equivalent to finding the companion matrix with multiple eigenvalues and nearest to $A=\begin{bmatrix}
        1&0\\
        1&0
    \end{bmatrix}$. We take $\mathcal{P}=\begin{bmatrix}
        e_1 & e_3
    \end{bmatrix} \in \C^{4 \times 2}$ and parametrize $u,v$ as in the previous example. Thus,
    \[   M = \begin{bmatrix}
        \sin(t) & e^{-i \theta} \cos(t)\\
    0&0\\
    \cos(t)&0\\
    0&\cos(t)
    \end{bmatrix}, \qquad r=\begin{bmatrix}
        (\lambda-1)\sin(t)\\
        \lambda e^{i \theta} \cos(t) - \sin(t)\\
        \lambda \cos(t) - \cos(t) + \sin(t) e^{-i \theta}\\
        -\lambda \sin(t) e^{-i\theta}
    \end{bmatrix}.\] Note first that, if $\cos(t)=0$, then necessarily $r \not\in \mathrm{Im}(M)$ and hence this choice is infeasible. If $\cos(t) \neq 0$, imposing $r \in \mathrm{Im}(M)$ yields the conditions $\theta \in \{0,\pi\}$ and $\lambda=\tan(t) e^{-i\theta}$. Under such conditions, we get
    \[ M^\dagger = \frac{1}{1+\cos(t)^2} \begin{bmatrix}
      \sin(t) & 0 & 2 \cos(t) & \sin(t) e^{-i \theta}\\
      \cos(t) e^{i \theta} & 0 & -\sin(t)e^{i \theta} & \frac{1}{\cos(t)}
    \end{bmatrix}, \qquad r=\begin{bmatrix}
        -\sin(t) [1 - e^{i\theta}\tan(t)] \\
        0\\
        2 \sin(t) e^{i \theta} - \cos(t)\\
        -2 \tan(t) \sin(t)
    \end{bmatrix}\]
    and we may indeed assume $\theta=0$ without loss of generality, as if $\theta=\pi$ we can just switch $t \leftarrow -t$. We ultimately get
    \[M^\dagger r = \begin{bmatrix}
        2 \tan(t)-1 \\
        -\tan(t)^2
    \end{bmatrix}.     \]
Setting $x=\tan(t)$ and differentiating $(M^\dagger r)^*M^\dagger r$ with respect to $x$, we easily see that $\| M^\dagger r\|$ is minimized by $t_0=\arctan(x_0)$  where $x_0 \approx 0.4534$ is the unique real solution to the cubic equation $x^3+2x-1=0$. We conclude that the companion matrix nearest to $A$ with a double eigenvalue (at $x_0$) is
\[  B = \begin{bmatrix}
    2 x_0 & -x_0^2\\
    1 & 0
\end{bmatrix}, \quad with \quad \ \|A-B\|_F = \sqrt{x_0^4 + (2 x_0-1)^2} \approx 0.2257,\]
and therefore the degree-two monic polynomial with a double root and nearest to $z^2-z$ is $z^2-2x_0  z +x_0^2=(z-x_0)^2$.
\end{example}

\section{Formulae in the unstructured case}\label{sec:2}
The formulae that we obtained in Section \ref{sec:problem} simplify significantly when $\mathcal{P} = I$ and $\delta=\vvec(\Delta)$, that is, $\mathcal{S}=\C^{n \times n}$ and no additional linear structure is imposed on the problem. In this case, we have that 
  \[ M=\begin{bmatrix}
  v^T \otimes I\\
  I \otimes u^*
  \end{bmatrix} \in\mathbb{C}^{2n\times n^2}. \]
It follows that
\begin{equation} \label{MMstar}
    MM^* + \varepsilon I = \begin{bmatrix}
    (1+\varepsilon)I& u v^T\\
    \overline{v}u^* &(1+\varepsilon)I
\end{bmatrix} =(1+\varepsilon)I + \begin{bmatrix}
    u&0\\
    0&\overline{v}
\end{bmatrix} \begin{bmatrix}
    0&v^T\\
    u^*&0
\end{bmatrix}.
\end{equation}
  Hence, by the Sherman-Morrison-Woodbury formula,
    \begin{equation} \label{eq:MMinv}
      (MM^* + \varepsilon I)^{-1} = \frac{1}{1+\varepsilon}I - 
      \begin{bmatrix}
    u&0\\
    0&\overline{v}
\end{bmatrix} 
\frac{1}{\varepsilon(1+\varepsilon)(2+\varepsilon)} \begin{bmatrix}1+\varepsilon & -1\\
-1 & 1+\varepsilon
\end{bmatrix}
\begin{bmatrix}
    0&v^T\\
    u^*&0
\end{bmatrix}.
    \end{equation}

Another useful observation is that $M(\overline{u} \otimes v)=0$ and $\begin{bmatrix}
    -u^* & v^T
\end{bmatrix}M=0$. This implies that $M$ is always singular, and hence in the problem without regularization ($\varepsilon=0$) the matrix $MM^*$ is not invertible. We can tell exactly what the eigenvalues of $MM^*$ are.
\begin{lemma}
    The eigenvalues of $MM^*$ are
    \[
    \sigma(MM^*) = (0,\underbrace{1,1,\dots,1}_{\text{multiplicity $2n-2$}},2).
    \]
\end{lemma}
\begin{proof}
    We know from~\eqref{MMstar} that $MM^*$ is a rank-2 modification of the identity, hence all except (at most) two eigenvalues are equal to $1$. There must be an eigenvalue $0$ since we have noted that $M$ has nontrivial kernel, and the last eigenvalue is determined by $\operatorname{Tr} MM^* = 2n$.
\end{proof}

Note that we have 
\[
\begin{bmatrix}
    -u^* & v^T
\end{bmatrix} r_1 = 0, \quad
\begin{bmatrix}
    -u^* & v^T
\end{bmatrix} r_0 = \varepsilon(v^T \overline{y_u} - u^*y_v),
\]
where we have set $y = \begin{bmatrix}
    y_v\\
    \overline{y_u}
\end{bmatrix}$. In particular, this implies that when $\varepsilon=0$ we have $r \in \operatorname{Im} M$, and hence the equation $M\delta = r$ is solvable. This observation ensures that the feasible region of~\eqref{fuvl} is non-empty for every $u,v,\lambda$.

\begin{corollary} \label{cor:abclambda}
In the notation above, define
\[
\alpha=r_0^*r_0, \quad
\beta=r_1^*r_0, \quad
\gamma_1 = v^T(A^T\overline{u}+\overline{y_u}\varepsilon), \quad
\gamma_2 = u^*(Av+y_v\varepsilon).
\]
Then, we have:
\[ a = \frac{2}{1+\varepsilon}, \qquad b=\frac{\beta}{1+\varepsilon}, \qquad c= \frac{\alpha}{1+\varepsilon} - \frac{\varepsilon(\overline{\gamma_1}\gamma_2 + \overline{\gamma_2}\gamma_1) - \abs{\gamma_1-\gamma_2}^2}{\varepsilon(1+\varepsilon)(2+\varepsilon)}. \]
\end{corollary}
\begin{proof}
Since $\begin{bmatrix}
    0&v^T\\
    u^*&0
\end{bmatrix}r_1=0$, multiplying~\eqref{eq:MMinv} by $r_1$ we obtain $(MM^*+\varepsilon I)^{-1} r_1 = \frac{1}{1+\varepsilon} r_1$.
Noting that $r_1^*r_1=2$, we compute
\[ a = r_1^* (MM^*+\varepsilon I)^{-1} r_1  = \frac{2}{1+\varepsilon};\]
\[ b = r_1^* (MM^*+\varepsilon I)^{-1} r_0  = \frac{\beta}{1+\varepsilon}.    \]
Finally, by multiplying $r_0$ to both sides of~\eqref{eq:MMinv}, we get
\begin{multline*}
c = r_0^*(MM^*+\varepsilon I)^{-1} r_0
= \frac{1}{1+\varepsilon}r_0^*r_0 -
\begin{bmatrix}
    \gamma_2\\
    \gamma_1
\end{bmatrix}^*
\frac{1}{\varepsilon(1+\varepsilon)(2+\varepsilon)} \begin{bmatrix}1+\varepsilon & -1\\
-1 & 1+\varepsilon
\end{bmatrix}
\begin{bmatrix}
    \gamma_1\\
    \gamma_2
\end{bmatrix}\\
= \frac{\alpha}{1+\varepsilon} - \frac{\varepsilon(\overline{\gamma_1}\gamma_2 + \overline{\gamma_2}\gamma_1) - \abs{\gamma_1-\gamma_2}^2}{\varepsilon(1+\varepsilon)(2+\varepsilon)}.
\end{multline*}
\end{proof}

The formula for $c$ here (and the resulting formula for the minimum $f_\varepsilon(u,v,\lambda_*) = \frac{ac-|b|^2}{a}$) is not particularly appealing, but we note that when $y=0$ it holds that $\gamma_1=\gamma_2$, and that $\lambda_*$ is independent of $\varepsilon$.

The minimizer also has an explicit form in the unstructured case.
\begin{proposition}  \label{prop:zuzv}
Fix $u,v,\lambda$ and let $\mathcal{P} = I$ (unstructured case), and $M,r$ be as in~\eqref{defM} and \eqref{defr}. 
Then, we have
\begin{equation} \label{zuzv}
    z_v = \frac{1}{1+\varepsilon}\left((\lambda I - A)v - y_v\varepsilon - u \eta_1\right), \quad
    z_u = \frac{1}{1+\varepsilon}\left((\lambda I - A)^* u - y_u\varepsilon - v\overline{\eta_2}\right).
`\end{equation}
where
\[
\eta_1 = \frac{(1+\varepsilon)\gamma_1-\gamma_2}{\varepsilon(2+\varepsilon)}, \quad \eta_2 = \frac{(1+\varepsilon)\gamma_2-\gamma_1}{\varepsilon(2+\varepsilon)},
\]
and $\Delta_* = z_v v^* + u z_u^*$.
\end{proposition}
\begin{proof}
Multiplying~\eqref{eq:MMinv} by $r = r_1\lambda + r_0$, we have
\[
z = \begin{bmatrix}
    z_v\\\overline{z_u}
\end{bmatrix} = (MM^*+\varepsilon I)^{-1}r = \frac{1}{1+\varepsilon} \begin{bmatrix}
    r_v\\\overline{r_u}
\end{bmatrix} 
+ 
\begin{bmatrix}
    u&0\\
    0&\overline{v}
\end{bmatrix} 
\frac{1}{\varepsilon(1+\varepsilon)(2+\varepsilon)}
\begin{bmatrix}1+\varepsilon & -1\\
-1 & 1+\varepsilon
\end{bmatrix}
\begin{bmatrix}
    \gamma_1\\
    \gamma_2
\end{bmatrix}.
\]
The formula for $\Delta_*$ then follows directly from Lemma~\ref{lem:minimizer_Prank2}.
\end{proof}

When $y=0$, the solution of the nearest matrix problem is a critical point for $f_{\varepsilon,0}$ irrespective of the value of $\varepsilon$. This is a special property of the specific nearness problem at hand, and it is potentially advantageous numerically; we prove it in Theorem \ref{thm:specialgradient}.
\begin{theorem} \label{thm:specialgradient}
    Let $u,v\in\mathbb{C}^n$ such that $u^* v = 0$ and $\norm{u}=\norm{v}=1$, and let $\delta>0$ be such that $(A-\lambda I) v = u\delta$ and $(A-\lambda I)^*u = v\delta$. Then, for every $\varepsilon > 0$, the point $[u,v]\in V_2$ is a zero of the Riemannian gradient of $f_{\varepsilon,0}$.
\end{theorem}
\begin{proof}
    We can compute using~\eqref{eq:MMinv} and~\eqref{defr}
    \[
    \begin{aligned}
        z &= (MM^*+\varepsilon I)^{-1}r\\
        &= \left(\frac{1}{1+\varepsilon}I - 
      \begin{bmatrix}
    u&0\\
    0&\overline{v}
\end{bmatrix} 
\frac{1}{\varepsilon(1+\varepsilon)(2+\varepsilon)} \begin{bmatrix}1+\varepsilon & -1\\
-1 & 1+\varepsilon
\end{bmatrix}
\begin{bmatrix}
    0&v^T\\
    u^*&0
\end{bmatrix}\right) \begin{bmatrix}
    -\delta u\\
    -\delta \overline{v}
\end{bmatrix}\\
&= -\frac{1}{1+\varepsilon}\begin{bmatrix}
    \delta u\\
    \delta\overline{v}
\end{bmatrix}
+
     \begin{bmatrix}
    u&0\\
    0&\overline{v}
\end{bmatrix} 
\frac{1}{\varepsilon(1+\varepsilon)(2+\varepsilon)} \begin{bmatrix}1+\varepsilon & -1\\
-1 & 1+\varepsilon
\end{bmatrix}
\begin{bmatrix}
    \delta\\
    \delta
\end{bmatrix}\\
&=\left(-\frac{1}{1+\varepsilon}+\frac{1}{(1+\varepsilon)(2+\varepsilon)}\right)
\begin{bmatrix}
    \delta u\\
    \delta\overline{v}
\end{bmatrix}\\
&=-\frac{1}{2+\varepsilon}
\begin{bmatrix}
    \delta u\\
    \delta\overline{v}
\end{bmatrix}
    \end{aligned}.
    \]
    Hence
    \[
    z_v = -\frac{\delta}{2+\varepsilon}u, \quad z_u = -\frac{\delta}{2+\varepsilon}v, \quad
    \Delta_* = z_vv^* + uz_u^* = -\frac{2}{2+\varepsilon}\delta uv^*.
    \]
    We now calculate the Euclidean gradient of $f_{\varepsilon,0}$ using~\eqref{gradf}; the first block is
    \[
    \begin{aligned}
    \nabla_u &= 2(\lambda_* I - A - \Delta_*)z_u\\
    &=\frac{2\delta}{2+\varepsilon} (-\delta u + \frac{2\delta}{2+\varepsilon}u)\\
    &= -\frac{2\delta^2\varepsilon}{(2+\varepsilon)^2}u
\end{aligned}
    \]
    and an analogous computation gives $\nabla_v = -\frac{2\delta^2\varepsilon}{(2+\varepsilon^2)}v$ for the second block.

    Importantly, the Euclidean gradient $g=[\nabla_u \ \ \nabla_v]$ is a multiple of $x = \begin{bmatrix}
        u & v
    \end{bmatrix}$, say, $g= \omega x$ for some $\omega \in \C$.
    The formula to compute the Riemannian gradient in the Stiefel manifold in a point $x\in\mathbb{C}^{n\times 2}$ starting from the Euclidean gradient $g\in\mathbb{C}^{n\times 2}$ is~\cite[Section~7.3]{Boumal}
    \[
    \operatorname{grad}_x f_{\varepsilon,0} = g - x H(x^*g),
    \]
    where $H(\cdot)$ denotes the Hermitian part. Plugging in $g = \omega x$, and recalling that $x^*x=I_2$, we get $\operatorname{grad}_x f_{\varepsilon,0} = 0$.
\end{proof}

\section{Numerical experiments}\label{sec:ne}

In this section, we test the performance of our algorithm in various settings. The implementation of the algorithms is available on~\url{https://github.com/fph/RiemannOracle}, in the functions named \texttt{nearest\_defective\_*}. The code relies on Manopt \cite{BoumalMishraAbsil} for performing optimization over the manifold $V_2$. For our experiments, we used  the function \texttt{penalty\_method} in the same repository as the solver, with the default options apart from the gradient tolerance \texttt{options.tolgradnorm = 1e-9}. The reason for this choice is that Theorem \ref{thm:specialgradient} suggests that the penalty method (i.e. setting $y=0$) might work especially well for this problem; nevertheless, our implementation also allows the user to call the augmented Lagrangian method, which is generally advertised as potentially better in the optimization literature.

\subsection{Choice of the starting point} \label{sec:startingpoint}

Unfortunately, the objective function $f$ has many local minima, and the one found by our algorithm depends on its starting point. This is analogous to what happens with the method from~\cite{ABBO}. Therefore, we adopt the same heuristic to select suitable starting points.

The heuristic relies on Theorem~\ref{thm:svdcharacterization}: namely, we choose a suitable candidate $\lambda_0$ for the multiple eigenvalue of $A+\Delta$, and then take $[u,v]$ to be the orthogonal factor in the thin QR factorization of $[u_n,v_n]$, where $u_n,v_n$ come from a singular value decomposition of $A-\lambda_0 I$. It only remains to select a suitable value of $\lambda_0$; for this we follow the heuristic strategy outlined in~\cite[page~506]{ABBO}: if $\lambda_1,\dots,\lambda_n$ are the eigenvalues of $A$ and $p_1,\dots,p_n$ are their condition numbers, then the pairs $\lambda_j, \lambda_k$ that are more likely to coalesce to a common eigenvalue of a perturbed matrix $A+\Delta$ are those with the smallest values of $s_{j,k} = \frac{\abs{\lambda_j-\lambda_k}}{p_j+p_k}$, and a first approximation of the coalescence eigenvalue is $\lambda_0 = \frac{p_j\lambda_k + p_k\lambda_j}{p_j+p_k}$. These heuristic formulas are justified by eigenvalue perturbation theory: A perturbation of small norm $\norm{\Delta}_F = \varepsilon$ moves $\lambda_j$ by a magnitude $\approx \varepsilon p_j$.

\subsection{Experiments with unstructured perturbations}
When applying the penalty method ($y=0$) to unstructured perturbations, we observe Theorem~\ref{thm:specialgradient} in practice: If, for a given value $\varepsilon=\varepsilon_0$ of the regularization parameter, the method converges to a critical point that satisfies the assumptions of Theorem~\ref{thm:specialgradient}, then the algorithm terminates immediately for all subsequent values of $\varepsilon \leq \varepsilon_0$. The method typically requires many iterations for the first value of $\varepsilon$, but then terminates immediately for later values.

We first benchmark our method on examples that exist in the literature for unstructured perturbations ($\mathcal{P}=I$); in particular, the article \cite{Russians23} provides some test problems with guaranteed exact solutions. We start from~\cite[Example~1]{Russians23}, i.e., the matrix
\[
A = \begin{bmatrix}
    1+i & 1-2i & 2-2i\\
    1+2i & 2+i & 1-3i\\
    2 & 1+2i & 2+i
\end{bmatrix}.
\]
With the starting point strategy outlined above, the pair with the lowest $s_{j,k}$ is $(j,k) = (1,2)$ and $\lambda_0 \approx -0.3393 + 1.2763i$; with this choice of starting point, the penalty method unfortunately converges to the local minimizer $\norm{\Delta_*}_F \approx 2.0886$. However, the global minimum corresponds to the coalescence of the eigenvalues $\lambda_2$ and $\lambda_3$ rather than $\lambda_1$ and $\lambda_2$. Selecting $(j,k)=(2,3)$ in the starting point heuristic (the choice with the second-smallest $s_{j,k}$) produces $\lambda_0 \approx 3.8109 + 0.6606i$, and converges to $\norm{\Delta}_F \approx 1.139495$, which is the global minimizer computed exactly via algebraic methods in~\cite{Russians23}. We note that the algorithm is not too sensitive to the choice of $\lambda_0$: in our experiments, any real choice of $\lambda_0 > 0.524$ gives the global minimum.

The picture is similar for~\cite[Example~2]{Russians23}, which is the flipped companion matrix
\[
A = \begin{bmatrix}
    0 & 1 & 0\\
    0 & 0 & 1\\
    -91 & -55 & -13
\end{bmatrix}.
\]
Choosing the pair $(j,k)$ with smallest $s_{j,k}$ gives $\lambda_0 = -3$, and the algorithm converges to a local minimum $\norm{\Delta_*}_F \approx 0.0836$. However, selecting once again a different pair of candidate coalescing eigenvalues $(j,k)=(2,3)$ produces the starting point $\lambda_0\approx -4.4680 - 1.2660i$ and the global minimum $\norm{\Delta_*}_F \approx 0.0350264$, which matches the results in~\cite{Russians23}.

The third and final example in~\cite{Russians23} is the $6\times 6$ Grcar matrix
\[
A = \begin{bmatrix}
    1 & 1 & 1 & 1 & 0 & 0\\ -1 & 1 & 1 & 1 & 1 & 0\\ 0 & -1 & 1 & 1 & 1 & 1\\ 0 & 0 & -1 & 1 & 1 & 1\\ 0 & 0 & 0 & -1 & 1 & 1\\ 0 & 0 & 0 & 0 & -1 & 1 
\end{bmatrix}
\]
(produced by Matlab's \texttt{gallery('grcar', 6}).
This exact example appears also in~\cite{AFS,AlamBora05,BGMN}. With the starting value $\lambda\approx 0.9850 - 0.5504i$ selected using the above heuristic (and the smallest $s_{j,k}$), we converge to the exact solution $\norm{\Delta_*}_F\approx 0.2151857666139$, which matches up to 12 significant digits the global minimum in~\cite{Russians23}.

Another family of experiments, taken from~\cite{AFS}, are the Kahan matrices produced by the Matlab command~\texttt{gallery('kahan', n, asin(0.1\^(1/(n-1))))}. For $n=6$, the starting values in~\cite{ABBO} produce only local minima, but selecting the same starting value $\lambda_0=0$ used in~\cite{AFS} gives the (likely global) minimum $\norm{\Delta_*}_F \approx 4.7049\times 10^{-4}$. However, for the Kahan matrix with $n=15$, both the  starting point from~\cite{ABBO} and the value $\lambda_0=0.12$ suggested in~\cite{AFS} produce only local minima (with the best value attained $\norm{\Delta_*}_F \approx 1.0031\times 10^{-6}$); in contrast, both the method in~\cite{AFS} and the code from~\cite{ABBO} find a better minimum $\norm{\Delta_*}_F \approx 4.4850\times 10^{-7}$. This confirms an observation already made in~\cite{RO}: The penalty method for spectral distance problems does not perform optimally when the minimum value $\min f(u,v)(=\norm{\Delta_*}_F^2)$ is close to machine precision. As possible explanations, we conjecture either numerical instability in our formulation of the objective function or practical issues in the Manopt implementation of the minimization subroutine. A further exploration of this aspect is an interesting direction for potential future improvement of our method. 

A larger-scale experiment, taken from~\cite{BGMN}, is the $67\times 67$ matrix \texttt{west0067} from Matrix Market~\cite{MatrixMarket}. Choosing the eigenvalue pair with the smallest $s_{j,k}$ produces the local minimum $\norm{\Delta_*}_F \approx 0.00602962$, while we need to take the fifth best $s_{j,k}$ to produce a starting point that gives convergence to the (likely global) minimum $\norm{\Delta_*}_F \approx 0.00551675$. These results coincide with the one from the code in~\cite{ABBO}, which also finds the local minimum $\norm{\Delta_*}_F \approx 0.00602962$ first, and then the global minimum $\norm{\Delta_*}_F \approx 0.00551675$ in its first two searches. For this $67\times 67$ matrix, computing the solution with the penalty method (for a single starting point $x_0$) takes about 2 seconds on a computer with Intel i5-1135G7 processor running Matlab R2024b.

In general, the new method performed comparably with state-of-the-art approaches on the selected experiments with dense unstructured perturbations. It did, however, tend to require more accurate starting points in order to capture the global minima, suggesting that it could particularly benefit from more reliable methods for generating initial guesses; we leave this aspect to future research. In the next subsection, we turn to structured perturbations, which represent the main novelty of this paper.

\subsection{Experiments with structured perturbations}

As a first simple test, we take the $2\times 2$  companion matrix in Example~\ref{example:companion}. The only eigenvalue pair yields $\lambda_0=0.5$, from which we obtain (up to numerical error) the exact solution computed in the example, with $\norm{\Delta_*}_F \approx 0.225713$.

We then take once again $A$ to be the Grcar matrix of size 6, and look for perturbations that preserve its Toeplitz structure (a subspace of dimension $p=11$). Using the eigenvalue pair with the two smallest $s_{i,j}$ to produce starting points, we obtain a local minimizer $\norm{\Delta_*}_F \approx 0.3180$. Using the next two smallest values of $s_{i,j}$, we obtain a better minimizer $\norm{\Delta_*}_F \approx 0.2309$; the corresponding matrix $A+\Delta_*$ is
\[
\begin{footnotesize}
\begin{bmatrix}
 1.0000 + 0.0001i & 1.0045 + 0.0060i & 1.0073 - 0.0027i & 0.9987 - 0.0034i & 0.0008 + 0.0016i & 0.0008 - 0.0036i\\
  -0.9799 + 0.0019i & 1.0000 + 0.0001i & 1.0045 + 0.0060i & 1.0073 - 0.0027i & 0.9987 - 0.0034i & 0.0008 + 0.0016i\\
 0.0108 + 0.0497i & -0.9799 + 0.0019i & 1.0000 + 0.0001i & 1.0045 + 0.0060i & 1.0073 - 0.0027i & 0.9987 - 0.0034i\\
  -0.0634 + 0.0481i & 0.0108 + 0.0497i & -0.9799 + 0.0019i & 1.0000 + 0.0001i & 1.0045 + 0.0060i & 1.0073 - 0.0027i\\
  -0.0811 - 0.0369i & -0.0634 + 0.0481i & 0.0108 + 0.0497i & -0.9799 + 0.0019i & 1.0000 + 0.0001i & 1.0045 + 0.0060i\\
  -0.0057 - 0.0735i & -0.0811 - 0.0369i & -0.0634 + 0.0481i & 0.0108 + 0.0497i & -0.9799 + 0.0019i & 1.0000 + 0.0001i\\
\end{bmatrix}
\end{footnotesize}
\]
with a double eigenvalue $\lambda_*\approx 0.7665 + 1.5825i
$.

We now take the Grcar matrix of size $n=15$, and look for perturbations that preserve both its Toeplitz and its zero structure; in this way, the subspace of allowed perturbations $\Delta$ has dimension $p=5$, since we can only change the numerical values for the five nonzero diagonals. This problem is more challenging than the previous tests, and for this reason we opted to choose different options for the solver, with a more cautious stepsize and a more conservative inner solver: \texttt{tolgradnorm=1e-8, epsilon\_decrease=0.8,solver=@trustregions}. With these options (and the initial value with the smallest $s_{i,j}$) we get $\norm{\Delta_*}_F \approx 0.2430$, and the values of the five nontrivial diagonals of $A+\Delta_*$ are $  -1.0071 - 0.0159i
, 1,    1.0030 - 0.0346i,    0.9570 + 0.0061i
,    1.0299 + 0.0152i$ respectively, starting from the subdiagonal. Note that it is expected that the main diagonal of $A+\Delta_*$ contains ones. Indeed, it is easy to prove that, if it contained a value $1+\alpha$ with $\alpha\neq 0$, then $A+\Delta_*-\alpha I$ would be another matrix with multiple eigenvalues, and nearer to $A$. It is therefore comforting that the algorithm returns such a solution. The computation requires about 3 seconds.

As a final test, we take once again the $67\times 67$ matrix \texttt{west0067}, which is sparse with 294 nonzeros, and we look for complex perturbations that preserve its sparsity structure; clearly, this is by construction a subspace of dimension $p=294$. We select the best starting point $\lambda_0 \approx   -0.2120 + 0.7296i$, and using the default options we obtain a sparse perturbation with $\norm{\Delta_*}_F \approx 0.0273$. The computation takes about 60 seconds, but for simplicity it was performed by allocating the full dense array $\mathcal{P}$ and without taking advantage of sparsity; we thus believe that, with a more careful implementation, it is possible to speed up the computation time.

\section{Conclusions}\label{sec:conclusions}
We have successfully modified the method in~\cite{RO} to compute the nearest matrix with multiple eigenvalues, a problem that has application to the stability theory of eigenvalues. Notably, the new method can deal with additional linear constraints on the perturbation matrix $\Delta$, such as imposing a sparsity or Toeplitz structure We integrated our code into the Matlab library \url{https://github.com/fph/RiemannOracle/}, and distributed it so that it is ready to use.

For unstructured perturbations, our new method performs comparably to state-of-the-art approaches, though it appears somewhat more dependent on the choice of initial points. However, our algorithm's ability to impose arbitrary linear constraints is a potentially valuable novelty,  with possible applications to the analysis of eigenvalue stability under structured perturbations. We believe that our research has also independent value in the development of the theory of Riemannian optimization methods for matrix nearness problems, because it extends the framework in~\cite{RO} with constraints on both left and right eigenvectors, rather than on right eigenvectors only.

 \section{Acknowledgements} We thank an anonymous reviewer and Michael Overton for sharing interesting comments and useful suggestions on the presentation.

\bibliographystyle{abbrvnat}
\bibliography{bibliography}

\end{document}